\documentclass{amsart}
\usepackage{hyperref}
\usepackage[]{color}
\usepackage{theorem}
\usepackage{calc}
\usepackage{ifthen}
\usepackage{epsfig}
\usepackage{graphicx}
\usepackage{latexsym}
\usepackage{amssymb,amsmath}

\newif\ifskip
\skipfalse
\newif\ifcut
\cuttrue
\newif\ifshort
\shortfalse

\begin{document}
\title{Modernism, Fiction and Mathematics}
\author{Johann A. Makowsky}
\maketitle
\begin{center}
\today
\end{center}
Review of:
\\
Nina Engelhardt,
\\
Modernism, Fiction and Mathematics,
\\
Edinburgh Critical Studies in Modernist Culture,
\\
Edinburgh University Press,
\\
Published June 2018 (Hardback), November 2019 (Paperback)
\\
ISBN
Paperback: 9781474454841, Hardback: 9781474416238

\section{The Book Under Review}
Nina Engelhardt's book is a study of four novels by three authors, 
Hermann Broch's trilogy {\em The Sleepwalkers} 
\ifshort
\else
\cite{broch1, broch1e},
\fi 
Robert Musil's {\em The Man without Qualities} 
\ifshort
\else
\cite{MoE,MoEe1,MoEe2},
\fi 
and Thomas Pynchon's  {\em Gravity's Rainbow} and {\em Against the Day} \cite{pynchon-GR,pynchon-AdD}.

Her choice of authors and their novels is motivated by the impact the mathematics of the interwar period had
on their writing fiction.
It is customary in the humanities to describe the cultural ambiance of the interwar period between World War I and World War II 
as modernism.
Hence the title of Engelhardt's book: Modernism, Mathematics and Fiction.
Broch and Musil are indeed modernist authors par excellence. Pynchon is a contemporary American author, usually classified by
the literary experts as postmodern. 

\ifshort
\else
\begin{figure}
\centering
\includegraphics[width=0.365\textwidth]{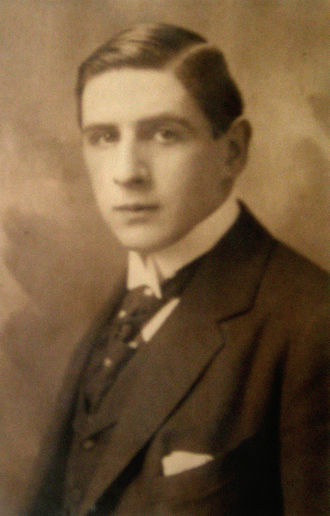} 
\hspace{2cm}
\includegraphics[width=0.4\textwidth]{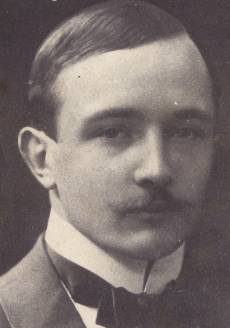}

Hermann Broch in 1909 \hspace{3.5cm}  Robert Musil in 1900
\end{figure}

\begin{figure}
\centering
\includegraphics[width=0.4\textwidth]{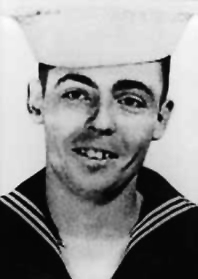}

Thomas Pynchon ca. 1957
\\
Works produced by employees of the United States federal government 
\\
in the scope of their employment are public domain by statute.

\end{figure}
\fi 


Engelhard's book is interesting for the literary minded mathematician for two reasons: First of all it draws attention
to three authors who spent a lot of time and thoughts in studying the mathematics of the interwar period and used this
experience in the shaping of their respective novels. Broch (1886-1951) and Musil (1880-1942)
 were Austro-Hungarians and wrote in German.
They were translated into English\footnote{The first translation of The Man Without Qualities in English was published by Ernst Kaiser and Eithne Wilkins in 1953, 1954 and 1960} 
immediately after World War II, 
and are recognized by connoisseurs and experts as the German language counterpart to Joyce and Proust.
However, even in the German language world they are not widely read,  
both because of their length of their novels and the demanding complex narrative.
Pynchon is also known for his lengthy novels and vast erudition, and is 
read mostly by his literary enthusiasts. 
All three authors deserve the attention of a wider audience among mathematicians and scientists,
as I will expound below.

The true merit of Engelhardt's book lies in its perspective. Without falling into the traps of over-interpretation, 
she discusses the interplay between modern mathematics and its foundational crisis and the development of the 
fictional characters described in the novels. 
She uses the modernist perspective cautiously to put the novels into the wider cultural context of
the prewar and interwar period  in central Europe, and she avoids successfully making wrong statements about 
the mathematics of this period.

If some readers of this review pick up only one of these novels and start reading it (possibly but not necessarily
together with the 
relevant chapters of Engelhardt's book) 
they will be amply rewarded. Hopefully these readers will then look for 
the necessary leisure time to go on reading more.
\section{Modernism}
\subsection{Modernity vs Modernism}
Modernism as a phenomenon in European culture is traditionally associated with the emancipation of art and science
from their exclusive role in representing the true view of the world based on the notion of immediate perception.
Kantian epistemology was questioned throughout the 19th century, a
cultural pessimism blossomed and the  church's monopoly on the sovereignty 
of interpretation of the true world was shattered.
Its breeding ground is the social changes in the period of intense industrialization before and after WWI. 

I do not think that there is a causal connection between the 
various manifestations of modernism. It is more natural and convenient to see it as a confluence of various 
developments of social, technological and even political character.

Certainly, some of these developments are interconnected, but not always in an obvious way.
The discovery of non-Euclidean geometry definitely was not a result of social or political changes,
neither was the emergence of Relativity Theory in physics, nor the discovery of the periodic table in chemistry,
nor the discovery of the Mendelian laws in biology.
On the other hand, discoveries in the natural sciences and technological progress are interrelated,
and there is no doubt that social and political developments were partially accelerated by technological advances.
\begin{itemize}
\item
In painting, 
the advent of photography, 
the rise of the  industrial bourgeoisie, the changing 
socio-economic conditions for painters, the confrontation with and the  influence of non-European art brought about
by colonialism, and the discovery
of non-Euclidean geometry, 
all may have played their role in the development of abstract painting.
Painting stopped to be the ``true picture'' of the world.
\item
In music, 
the advent of acoustic recording played the same role as photography in painting, and the new discoveries of
Helmholtz and others in acoustics played a similar role as the discovery of 
non-Euclidean geometry. 
The encounter with non-European music together with the changing social standing of professional musicians
opened the ways for experimentation. Twelve-tone music was not only developed by Arnold Schoenberg and 
his followers in Vienna 
but also by Nikolai Roslavets in Russia. 
\item
In literature the modern and modernist trends are more varied and it is  difficult to bring them into a 
coherent picture. 
Trends in poetry, drama, and prose are diverse.
Influences of philosophy are divergent. Schopenhauer, Nietzsche, Bergson and Kierkegard were among the
influential cultural pessimists. Ernst Mach, Carl Stumpf, Edmund Husserl, Ernst Cassirer, the Neo-Kantians, and the
philosophers of the  Vienna Circle were reexamining the nature of perception and the scientific method, and 
altered our epistemological outlook. 

\ifskip
\else
Peter Szondi's {\em Theorie des modernen Dramas} \cite{szondi}
and Martin Esslin's {\em The Theater of the Absurd} \cite{esslin},
and George Lukacs's
{\em The Theory of the Novel}
\ifshort
\cite{lukacs1}
\else
\cite{lukacs1,lukacs2}
\fi 
and Milan Kundera's {\em The Art of the Novel} \cite{kundera},
may be viewed as landmarks in the theoretical understanding of modern drama
and modern novels.
\fi 

For our discussion here it is noteworthy that among the German language authors of
novels  
Hermann Broch,
Hermann Hesse, Thomas Mann,  and Robert Musil (in alphabetical order)
are considered paradigmatic for modernist prose in German. The same holds for James Joyce in English 
and Marcel Proust and Paul Val\'ery in French.
\item
Modern mathematics is  largely the result of two internal developments in mathematics: 
the discovery of non-Euclidean geometry and the emergence of set theory.
In the years between 1874 and 1897 Georg Cantor created set theory single-handed, followed by Felix Hausdorff's
{\em Grundz\"uge der Mengenlehre} \cite{hausdorff}
and
in 1899 Hilbert published his {\em Foundations of Geometry} 
\cite{hilbert} and started emphasizing the use of 
\ifshort
\else
axiomatic methods, see \cite{corry2004david}.
\fi 
To the readers of the {\em NOTICES},
``modernism in mathematics'' may sound familiar from Yuri Manin's review \cite{manin}
of Jeremy Gray's {\em Plato's Ghost: The Modernist Transformation of Mathematics} \cite{gray},
and from  Frank Quinn's
{\em Invisible Revolution of Modern Mathematics} \cite{quinn}.
In his review of \cite{gray} Solomon Feferman \cite{feferman} describes ``modern mathematics'' as follows:
\begin{quote}
\em
Modern mathematics--in the sense the term is used by working mathematicians these days--took shape 
in the period from 1890 to 1930, mainly in Germany and France. 
Strikingly new concepts were introduced, new methods were employed, 
and whole new areas of specialization emerged, while other themes were relegated to the dusty shelves of history.
At the same time, the nature of mathematical truth and even the consistency of mathematics were put into question,
as mathematicians, logicians and philosophers grappled with the subject's very foundations.
[...] modernism is defined as an autonomous body of ideas, having little or no outward reference, 
placing considerable emphasis on formal aspects of the work and maintaining a complicated--indeed, 
anxious--rather than naive relationship with the day-to-day world.
\end{quote}
Books reflecting the spirit of modern mathematics are sometimes called 
``Modern Algebra'' (van der Waerden), ``Modern Analysis'' (Dieudonn\'e, Rudin)
and the modern approach culminates in Bourbaki's unfinished multi-volume enterprise.
%
%
\item
In the realm of classical Number Theory
one might be tempted to consider the Prime Number Theorem of Hadamard and de La Vall\'ee-Poussin
as the first success of modern mathematics for its use of analytic functions in the proof.
Hardy's belief that no  ``elementary proof'' of it should exists reflects his belief 
that modern methods give rise to deeper mathematical understanding and can not be dispensed of. 
H\'elas, Erd\H{o}s and Selberg shattered this belief. For a history of the Prime Number Theorem see \cite{Goldstein}.
Similarly, Wiles' proof of the Shimura conjectures is just one of the later achievements of modern
mathematics, leading to a proof of Fermat's Last Theorem.  Its original proof uses Grothendieck's language
of schemes and universes.
To what extent this can be eliminated is discussed by A. MacIntyre \cite{macintyre}
and McLarty \cite{mclarty}.
The question to what extent Fermat's Last Theorem
admits an  ``elementary proof'' is still open.
\end{itemize}
In order to appreciate Nina Engelhardt's book we can summarize the modernist aspect of (modern) mathematics
as follows: The intuitive notion of mathematical truth (facts) is replaced by the axiomatic description
of relationships where the only concern is consistency. Mathematics passes from describing the true world
to describing mere possibilities.

\section{Formalists vs intuitionists}
One of the first modern(ist) mathematicians who was also a successful philosopher, playwright and essayist
was Felix Hausdorff.
He started his career as a successful writer under the name of Paul Mongr\'e, and his writings could easily
be classified as modernist. In 1903/04, in his lecture course {\em Time and Space} he remarked 
about mathematics in general,
reflecting on Hilbert's view on the axiomatic method:
\begin{quote}
{\em
Mathematics stands completely apart not only from the actual meaning that one attributes to its concepts but also from the 
actual validity one ascribes to its propositions. Its undefinable concepts are arbitrarily chosen objects of thought, its 
axioms are also arbitrary, though chosen so as to be free from contradiction. 
Mathematics is the science of pure thought, just as formal logic.}
\footnote{
Quoted from \cite{purkert}
}
\end{quote}

His literary publications were mostly published between 1897 and 1904 and consisted of a volume of poems, 
a play, a book on epistemology,
and a volume of aphorisms. He wrote no novel, and as a philosopher he was mostly influenced by the younger Friedrich Nietzsche. 
However, he was in correspondence with Moritz Schlick, the founder of the Vienna Circle,
and seemingly influenced Schlick's early views on time and space.

The modern(ist) view of mathematics
led to heated debates in continental Europe between Formalists and Intuitionists, with heavy political
overtones.
In philosophical circles such as the Vienna Circle around Hans Hahn and 
Moritz Schlick, and its Berlin version around Hans Reichenbach this view led to a new and very 
strict definition of what constitutes the scientific method.
The two philosophical circles were also influenced by
Ludwig Wittgenstein and the emergence of mathematical logic.
Both Kurt G\"odel and Karl Menger were students of Hans Hahn.
Among the eminent mathematicians
Karl Menger was an active member of the Vienna Circle and David Hilbert became a member of the Berlin Circle.

\section{Modernist writers and mathematics}
While I studied mathematics in the late 1960s 
I was preoccupied by the question of the nature of mathematical thought. I read a lot of material pertaining to this question:
philosophy, psychology, neurology, and also literary works.
Among the latter
I came across five authors who struck me as significantly
influenced by modern mathematical and/or musical thought, Paul Val\'ery, Robert Musil, Hermann Broch and Hermann Hesse,
and to a lesser extent Thomas Mann.
All of them are Europeans who wrote mostly before WWII.
Broch and Musil are the subjects of Engelhardt's book. 
Val\'ery, Mann and Hesse are not. 
\ifshort
\else
I will briefly discuss, whether they should have been.
\fi 
Thomas Pynchon is American and contemporary and writes in English. Engelhardt's book originated as a Ph.D. thesis from 2012
in the department of English Literature at the University of Edinburgh, 
{\em Mathematics in literature: modernist interrelations in novels by Thomas Pynchon, Hermann Broch, and Robert Musil}. 
Pynchon's two novels both deal with the mathematics
and science, and are set partially  in continental Europe of the period of Broch and Musil.
I will also briefly discuss why Pynchon, besides fitting a thesis in the department of English Literature, is
a major author discussed in Engelhardt's book.

\ifskip
\else
\subsection*{
Paul Val\'ery, poet and philosopher} 
Paul Val\'ery,
a symbolist/modernist French man of letters is best known for his poetry. However, he was also
interested in philosophy and mathematics, and was supposed to write a preface for F. Le Lionnais' {\em Les grands courants de la 
pens\'ee math\'ematique}, but died in 1945 before writing his preface. 
Raymond Poincar\'e, Louis de Broglie, Andr\'e Gide, Henri Bergson, and Albert Einstein all appreciated Val\'ery's 
thinking. They became friends personally and by correspondence. 
His {\em L'id\'ee fixe ou deux hommes \`a la mer} 
\ifshort
\cite{valery}, 
\else
\cite{valery,valery-e}, 
\fi 
written in 1932, 
is a brilliant reflection on the relationship between art and science,
where he describes the gap between experience and theory as so big that 
only a visionary architect can possible bridge it.
The text resonates Val\'ery's  experience from attending a lecture by 
Albert Einstein on Time, Space and Relativity.
\ifshort
\begin{figure}
\centering
\includegraphics[width=0.2\textwidth]{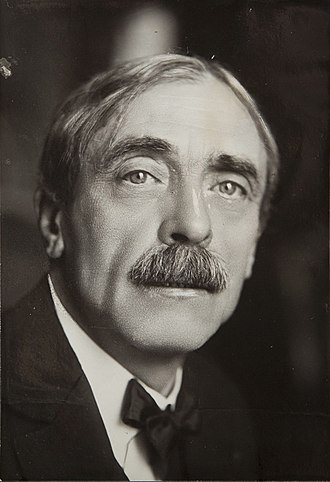}

Paul Val\'ery \\ picture taken by Henri Manuel, ca. 1925.
\end{figure}
\else
\begin{figure}
\centering
\includegraphics[width=0.4\textwidth]{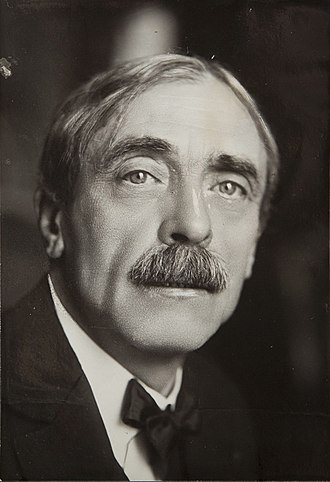}

Paul Val\'ery, picture taken by Henri Manuel, ca. 1925.
\end{figure}
\fi 
\fi 

\subsection*{Hermann Broch and the Vienna Circle} 
Hermann Broch was born in 1886 in Vienna
to a prosperous Jewish family and worked for some time in his family's factory.
As the oldest son, he was expected to take over his father’s textile factory, therefore, 
he attended a technical college for textile manufacture and a spinning and weaving college.
Although he had from early on philosophical, mathematical and literary interests, he pursued them only privately,
mostly at night, while managing the textile factury during the day.
In the years 1904-1905 he attended seminars and public lectures in mathematics and philosophy at Vienna University,
especially those given by Ludwig Boltzmann.

After an unsuccessful marriage which led to his divorce in 1923, he sold the textile factory in 1927 to devote himself
full time to his intellectual interests. He studied mathematics, philosophy, physics and also psychology and cultural history
at the University of Vienna. Philosophically he was, like Hausdorff, first under the influence of Schopenhauer, Nietzsche
and Weininger, but later turned to neo-Kantian philosophy and logicism and got involved with the Vienna Circle,
siding with Wittgenstein and Menger rather than being a strict positivist. From 1913 on he started to publish
philosophical essays. 

\ifskip
\else
Two recent theses which appeared during or after Engelhardt's work on her book are devoted in great
detail to Broch's philosophical and scientific background and its impact on his literary work:
{\em Ludwig Wittgenstein and Hermann Broch: The Need for Fiction and Logic in Moral Philosophy}
by
C. W. Bailes \cite{bailes}
and {\em Zweige eines einzigen Stamms. Hermann Broch's Reflections on Science and Literature} by
H.A.M. Brookhuis \cite{brookhuis}.
\fi 
Broch studied at Vienna University until 1930 having among his teachers 
Moritz Schlick, Rudolf Carnap, Hans Hahn and Karl Menger.
\ifshort
\else
\begin{figure}
\centering
\includegraphics[width=0.18\textwidth]{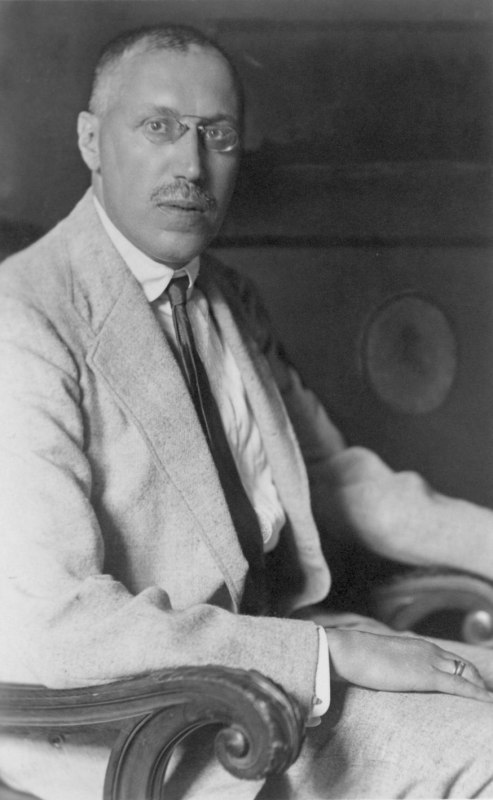}
\hspace{2cm}
\includegraphics[width=0.22\textwidth]{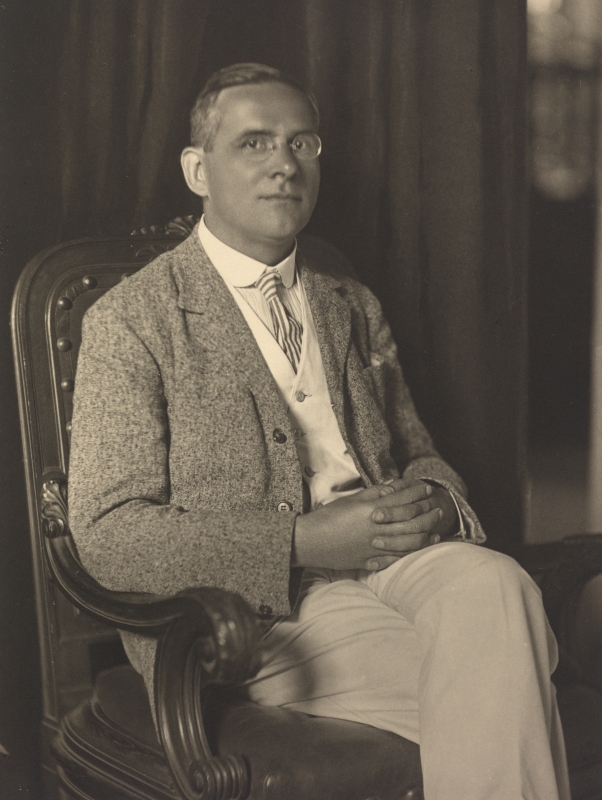}
\hspace{2cm}
\includegraphics[width=0.21\textwidth]{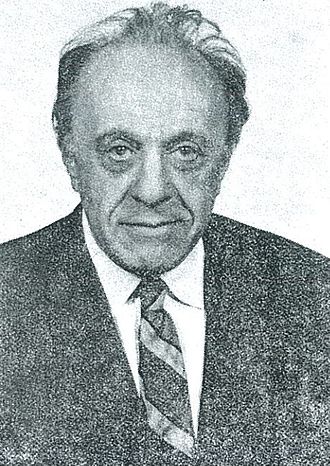}

\hspace{0.2cm} 
Hans Hahn 
\hspace{2.5cm} 
Moritz Schlick \hspace{2.8cm} Karl Menger
\end{figure}
\fi 
Originally he was more interested in a scientific career and wanted to write a Ph.D., 
but he became discouraged, due his lack of knowledge in Latin, which was, at that time, 
a compulsory subject for getting a Ph.D.
The goal of writing a thesis in mathematics was delayed indefinitely and 
finally abandoned also because of his disenchantment with science due to its
increased specialization.
It is during this period that he wrote his trilogy {\em The Sleepwalkers} published in 1931 which marked the beginning
of his career as a novelist, and which became known as a modernist novel par excellence.
He also wrote a novel, {\em Die Unbekannte Gr\"osse} 
\ifshort
\else
\cite{broch2,broch2e}
\fi 
inspired by the events around Einstein's theory and its experimental verification during an eclipse visible only
in the southern hemisphere.
Among the writers which N. Engelhardt discusses, H. Broch is definitely the most educated mathematician.

\subsection*{Robert Musil and the Berlin Circle}
Robert Musil  was born in 1880 in Klagenfurt (Austria) into a middle class family. His father was a professor of 
mechanical engineering at the German Technical University in Br\"unn (now Brno, Czech Republic), and, later, 
rose to hereditary nobility in the Austro-Hungarian Empire.
Musil held two doctorates, one in mechanical engineering at his father's university (1901), and one in
in psychology and philosophy at the University of Berlin under the renowned Professor Carl Stumpf.
Carl Stumpf wrote his Habilitationsschrift in 1870 in G\"ottingen on the foundations of mathematics.
\ifshort
\else
It was
recently republished  \cite{stumpf}.
\fi 

\ifshort
\else
\begin{figure}
\centering
\includegraphics[width=0.4\textwidth]{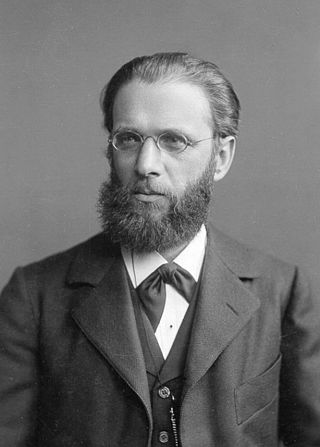}

Carl Stumpf, ca. 1900
\end{figure}
\fi 
M. Schlick's early philosophy shows both appreciation and criticism
of Stumpf's epistemology.
Among Stumpf's students we also find Edmund Husserl, who established the school of phenomenology,
and who had studied mathematics with Kronecker and Weierstrass.
Musil's majors in Berlin were ``logics and experimental psychology'', and he was examined in 1908 in
philosophy, natural science and mathematics.
A. Schwarz of the renowned Schwarz' inequality was among his doctoral examiners.
The topic of his thesis was {\em Beitr\"age zur Beurteilung der Lehren Machs}.
His thesis and other technical writings, covering the period from 1904-1922 were republished in 1980 as \cite{musil-phd}.
In 1906 his first novel was published and at the same time  he was developing an apparatus 
to research how people experience color.
From then on he dedicated his time exclusively to his literary career, 
interrupted only by his military service during WWI.
From 1918 till 1931 Musil lives in Vienna.
We do not know whether Musil had direct contact with the Vienna Circle during these years in Vienna.
In 1931-33 Musil lived in Berlin and was an active participant in the Berlin Circle around Hans Reichenbach and Richard von Mises. 
However, the first two volumes of {\em The Man Without Qualities} were already published in 1930.
Musil's awareness of modern mathematics and the foundations of geometry must have been acquired while he studied for 
his doctorate in Berlin. 
Musil's background in engineering, psychology and mathematics most definitely influenced his writing not only in 
{\em The Man without Qualities}.
In his 1906 novel {\em The Confusions of Young T\"orless} the confusions are manifold, 
ranging from the sexual to the foundations of geometry.

\subsection*{Thomas Pynchon's failed attempt to get a Ph.D. in mathematics}
Thomas Pynchon is famously reclusive.
About his mathematical and philosophical education we know rather little.
As a writer he belongs more to post-modernism than to modernism.

I gather from the Wikipedia article about him the following:
He was born in 1937 and holds an English degree from Cornell.
Starting in 1953 he
studied engineering physics at Cornell University, but left at the end of his second year to serve in the U.S. Navy. 
In 1957, he returned to Cornell to pursue a degree in English.
From February 1960 to September 1962, he was employed as a technical writer at Boeing in Seattle. 
In 1964, his application to study mathematics as a graduate student at the University of California, Berkeley was turned down.
The two novels discussed in Engelhardt's book are
{\em Gravity's Rainbow}, published in 1963, and {\em Against the Day}, published in 2006.
The former draws from his experience at Boeing, the latter contains a lot of mathematical material pertaining
to Sofia Kovalevskaya
and to Hilbert's school in G\"ottingen. Pynchon seemingly researched this material with the help of Michael Naumann, 
German secretary of culture from 1998 until 2001.
The only other good source for Pynchon's biographical material is his lengthy autobiographical introduction
to the collection of short stories {\em Slow Learner} \cite{pynchon-SL}, published in 1984.

The novels of Pynchon are considered difficult to read. Nevertheless, he has his fans who created a
Wikipedia-like website ({\tt https://pynchonwiki.com/}), which should help potential readers to
enter Pynchon's world. Literary critiques celebrate him as one of the most important American novelists of
the second half of the 20th century. {\em Gravity's Rainbow} and, even more so, {\em Against the Day}
are novels fitting the topic of Engelhardt's study. Contrasting him to Broch and Musil, 
comparing modernist to post-modernist writing, makes Engelhardt's book a gem among the few studies linking
mathematics and literature.
\ifskip
\else
\subsection*{Hermann Hesse and Thomas Mann}
There are two German authors who come to my mind for their modernist novels 
dealing with formal aspects of mathematics and music:
Thomas Mann with his {\em Doktor Faustus}
and Hermann Hesse with his {\em Magister Ludi (aka The Glass Bead Game)}.
Thomas Mann's {\em Doktor Faustus} does not deal with Mathematics and its foundational crisis but instead with music
and the atonality of the Vienna School. Thomas Mann was not an expert in modern music. 
He commissioned
Theodor Adorno to help with the musical aspects of his novel. But in a a very precise sense his novel could have been discussed
in the Framework of Modernism, Music and Fiction. It was wise of Nina Engelhardt not to include Thomas Mann's novel 
in her discussion, but it may be a tempting topic for further research.
\fi 

\ifshort
\else
\begin{figure}
\centering
\includegraphics[width=0.4\textwidth]{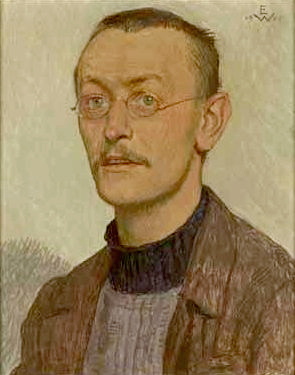}

Hermann Hesse (1877-1962) 

1905 portrait by Ernst W\"urtenberger (1868-1934)
\end{figure}
\fi 

\ifskip
\else
In Hermann Hesse's novel the situation is more complex. 
\fi 
One of the central themes of the novel is 
the {\em Glass Bead Game} and the way it is performed in public by its masters.
The details of the game itself are only vaguely described in the novel as an 
amalgam of mathematics, music and the Chinese I Ching.
But it says that
playing the game well requires years of hard study of music, mathematics, and cultural history. 
The game is essentially an abstract synthesis of all arts and sciences. 
It proceeds by players making deep, and sometimes unexpected, connections between seemingly unrelated topics.
The novel is a dystopia in a post-historic setting.
The mathematics is not specified in detail, but the sources used by Hesse concerning music,  history, and Indian and Chinese
culture are well documented and identified. 

Carl Jacob Burckhardt, the Swiss historian and diplomat, 
of the same family as  the eminent historian and 
scholar of Renaissance culture Jacob Burckhardt, 
was a  close friend of Hesse and can be recognized as one of the novel's protagonists. 
His career alternated between periods of academic historical research and diplomatic postings; 
the most prominent of the latter were League of Nations High Commissioner for the Free City of Danzig (1937–39) 
and President of the International Committee of the Red Cross (1945--48). [...] 
He was teaching contemporary history
at the University of Zurich from 1927 until 1932\footnote{
Quoted from the wikipedia {\tt https://en.wikipedia.org/wiki/Carl\_Jacob\_Burckhardt}}.
\ifshort
\else
\begin{figure}
\centering
\includegraphics[width=0.35\textwidth]{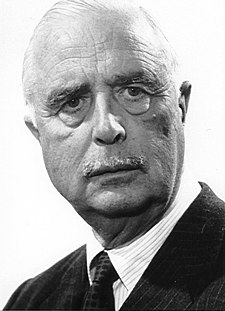}

Carl Jacob Burckhardt (1891-1974), 
\\
former President of the International Committee 
\\
of the Red Cross (ICRC) from 1944 to 1948
\\
(copyright ICRC)
\end{figure}
\fi 

In the Hesse archive there is no
written evidence of Hesse's inspiration concerning modern mathematics.
However, in a talk given by Hermann Weyl,
\cite{weyl-zurich} we learn, that Hermann Weyl, Hermann Hesse and Carl Jacob Burckhardt were close friends.

\ifshort
\else
\begin{figure}
\centering
\includegraphics[width=0.4\textwidth]{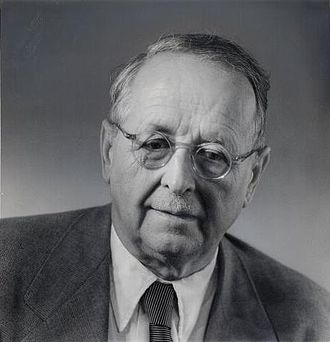}

Hermann Weyl, 1885--1955
\\
copyright ETH-Bibliothek Z\"urich, Bildarchiv
\end{figure}
\fi 

Commenting on his decision to leave Zurich to succeed Hilbert in G\"ottingen in 1930, Weyl writes:
\begin{quote}
\em
[...] it could have possibly been right to stay in Zurich, where my position was specially 
favorable to a vita contemplativa, paired with, relatively speaking, modest external impact.
It was not easy for me to come to a clear decision. I have pondered the question in my mind together with two
men, with Jacob Burckhardt and Hermann Hesse. The fate of the latter is particularly close to me. He went ``South'', not only
geographically. His example is tempting, but also a warning. 
For he achieved only a deeper loneliness and a mild desperation, which strangely permeated the happiness of his senses, 
his eyes and his writing, which was at the same time exhilarating and loosening up. 
Burckhardt on the other hand rejected the offer to succeed Ranke in Berlin with self-evident determination.
\end{quote}
Hermann Weyl was known for his inspiring lectures, celebrating mathematics as a performing art.
It is then likely, that he was both the model for the Magister Ludi celebrating the glass bead game,
and the source for Hesse's awareness of modern mathematics as the art of the possible, rather than the science of true nature.
The impact of Hermann Weyl on Hesse's Magister Ludi has not been explored in depth, but many mathematicians reading the novel
conjecture modern formalist mathematics to be a major source for Hesse's concept of the Glass Bead Game.
Few know that the source of Hesse's knowledge seems to be Hermann Weyl.

\section{What Engelhardt achieves}
\ifshort
\else
The book under review has no abstract, but the underlying Ph.D. thesis \cite{engelhardt-phd}, 
from 2012 does. It says:
\begin{quote}
\em
The focus of this thesis is on four novels’ illustrations of the parallels and interrelations between the 
foundational crisis of mathematics and the political, linguistic, and epistemological crises around 
the turn to the twentieth century.
While the latter crises with their climax in the First World War 
are commonly agreed to define modern culture and literature, this thesis concentrates on their relations 
with the ``modernist transformation'' of mathematics as illustrated in 
Thomas Pynchon’s Against the Day (2006) and Gravity’s Rainbow (1973), 
Hermann Broch’s The Sleepwalkers (1930-1932), and 
Robert Musil’s The Man without Qualities (1930/32). 
\end{quote}
\fi 
I have explained why the choice of these four novels is justified by examining the mathematical background
of the three authors Broch, Musil and Pynchon.
I have also indicated that, in a further study along the same lines,
a close examination of Hesse's {\em Glass Bead Game} could be a valuable addition.

Engelhardt's view of the ``modernist transformation'' of mathematics is close to what Quinn describes in \cite{quinn},
although slightly less precise. In her words  the impact of this ``modernist transformation'' on fiction consists in the
development of autonomy and independence from the natural world, almost a kind of escapism into possible worlds.
Musil's {\em man without qualities} is a {\em man of all possibilities}.

\ifshort
\else
In the abstract published by Edinburgh University Press advertising the book on their website, it says further:
\begin{quote}
\em
[The four novels] accord mathematics and its modernist transformation a central place 
in their visions and present it as interrelated with political, linguistic, 
epistemological and ethical developments in the modern West. 
Not least, the texts explore the freedoms and opportunities that the mathematical crisis implies 
and relate the emerging notion of ‘fictional’ characteristics of mathematics to the possibilities of literature. 
By exploring how the novels accord mathematics a central role as a particularly 
telling indicator of modernist transformations, 
this book argues that imaginative works contribute to establishing mathematics as part of modernist culture. 
\end{quote}
\fi 
Engelhardt's
monograph offers the readers a new framework 
of textual and cultural analysis in which one can
understand the modernist and postmodern
interplay
of literature and mathematics.
It is fair to say that in her book Nina Engelhardt does succeed in
giving us an inspiring tour d'horizon
of this interplay.

\bibliographystyle{amsplain}      
\bibliography{a-ref}   
\newpage
\appendix
\small
\section*{
Souces and copyright of the illustrations}
In order of their appearence.
\begin{itemize}
\item[Hermann Broch:]
public domain
\\
https://en.wikipedia.org/wiki/Hermann\_Broch
\item[Robert Musil:]
public domain
\\
https://en.wikipedia.org/wiki/Robert\_Musil
\item[Thomas Pynchon:]
public domain
\\
https://www.vulture.com/2013/08/thomas-pynchon-bleeding-edge.html
\\
Credit:
Pynchon, age 16, in his 1953 high-school yearbook, one of the few known photos of the author. 
Photo: Getty
\\
https://www.cs.mcgill.ca/~rwest/wikispeedia/wpcd/images/281/28137.jpg.htm
\item[Paul Valery:]
public domain
\\
https://en.wikipedia.org/wiki/Paul\_Val\'ery
\item[Hans Hahn:]
public domain
\\
https://geschichte.univie.ac.at/en/persons/hans-hahn-prof-dr
\\
Hans Hahn (1879-1934), Mathematics
Courtesy: Archive of the University of Vienna, picture archive Originator: Theo Bauer, Wien XIV. Signatur: 106.I.1210
1927
\item[Moritz Schlick:]
public domain
\\
https://en.wikipedia.org/wiki/Moritz\_Schlick
\item[Karl Menger:]
public Domain
\\
https://en.wikipedia.org/wiki/Karl\_Menger
\\
Photograph of mathematician Karl Menger, from the January 1970 issue of the Shimer College Bulletin.
Shimer College - 
\\
https://archive.org/details/ShimerCollegeBulletin\_1970-01 , page 1
\item[Carl Stumpf:]
public domain
\\
https://en.wikipedia.org/wiki/Carl\_Stumpf
\item[Hermann Hesse:]
public domain
\\
https://en.wikipedia.org/wiki/Hermann\_Hesse\#/media/File:Hermann\_Hesse.jpg
\\
1905 portrait by Ernst Würtenberger (1868–1934)
\item[Carl Jacob Burckhardt:]
public domain
\\
https://en.wikipedia.org/wiki/Carl\_Jacob\_Burckhardt
\\
Copyright: http://www.redcross.int/EN/copyright.html
\\
Any part of this website can be cited, copied, translated into other languages or adapted to meet 
local needs without prior permission from any component of the Movement, provided that the source is clearly stated.
\item[Hermann Weyl:]
public domain
\\
https://en.wikipedia.org/wiki/Hermann\_Weyl
\\
Source	ETH-Bibliothek Zürich, Bildarchiv, CC BY-SA 3.0
\\
https://en.wikipedia.org/wiki/Hermann\_Weyl\#/media/File:Hermann\_Weyl\_ETH-Bib\_Portr\_00890.jpg
\end{itemize}

\end{document}